\documentclass[12]{article}
\DeclareMathAlphabet{\mathpzc}{OT1}{pzc}{m}{it}
\usepackage{amsmath}
\usepackage{graphicx}
\begin{document}

\Large
\begin{center}
\bf{Doily as Subgeometry of a Set of Nonunimodular\\ Free Cyclic Submodules}
\end{center}
\vspace*{-.3cm}

\large
\begin{center}
 Metod Saniga$^{1}$ and Edyta Bartnicka$^{2}$
\end{center}
\vspace*{-.5cm} 

\normalsize
\begin{center}

$^{1}$Astronomical Institute of the Slovak Academy of Sciences,\\
SK-05960 Tatransk\' a Lomnica, Slovak Republic\\
(msaniga@astro.sk)

and

$^{2}$University of Warmia and Mazury,
Faculty of Mathematics and Computer Science, S\l{}oneczna 54 Street,
P-10710 Olsztyn, Poland\\ (edytabartnicka@wp.pl)

\end{center}

\vspace*{-.4cm} \noindent \hrulefill

\vspace*{-.1cm} \noindent {\bf Abstract}

\noindent
It is shown that there exists a particular associative ring with unity of order 16 such that the relations between nonunimodular free cyclic submodules of its two-dimensional free left module can be expressed in terms of the structure of the generalized quadrangle of order two.  Such a doily-centered geometric structure is surmised to be  of relevance for quantum information.\\

\vspace*{-.2cm}
\noindent
{\bf Keywords:}  associative ring with unity -- free cyclic submodules -- generalized quadrangle  


\vspace*{-.2cm} \noindent \hrulefill

\bigskip
\noindent
Let $R$ be a finite associative ring with unity ($1$) and $R^2$ its free left module. The set
$R(a,b) = \left\{ (\alpha a, \alpha b) | (a,b) \in R^2, \alpha \in R   \right\}$
is called a cyclic submodule of $R^2$. $R(a,b)$ is called {\it free} if the mapping $\alpha \mapsto (\alpha a, \alpha b)$ is injective.
A pair/vector ($a, b$) $\in R^{2}$ is called {\it unimodular} over $R$ if there exist $c, d \in R$ such that $ac + bd = 1$. 
It is well known (see, for example, \cite{veld}) that if $(a,b)$ is unimodular, then $R(a,b)$ is free. A great majority of finite rings have the property that all their free $R(a,b)$'s are generated by unimodular vectors. Here we shall consider a specific ring where this is not true, that is, a ring that also features free $R(a,b)$'s containing no unimodular vector; in what follows we shall call such free cyclic submodules nonunimodular. 
Our ring $R$ is a non-commutative one of order 16, defined as follows:  

\begin{equation*}
    R = \left\{ \left(
\begin{array}{ccc}
a & c & d \\
0 & b & 0 \\
0 & 0 & b \\
\end{array}
\right) \mid ~ a, b, c, d \in GF(2) \right\}.
\end{equation*}
Labeling the 16 matrices as follows
\begin{equation*}
~~0 \equiv \left(
\begin{array}{ccc}
0 & 0 & 0 \\
0 & 0 & 0 \\
0 & 0 & 0 \\
\end{array}
\right),~~
1 \equiv \left(
\begin{array}{ccc}
1 & 0 & 0 \\
0 & 1 & 0 \\
0 & 0 & 1 \\
\end{array}
\right),~~
2 \equiv \left(
\begin{array}{ccc}
1 & 1 & 0 \\
0 & 1 & 0 \\
0 & 0 & 1 \\
\end{array}
\right),~~
3 \equiv \left(
\begin{array}{ccc}
0 & 1 & 0 \\
0 & 0 & 0 \\
0 & 0 & 0 \\
\end{array}
\right),
\end{equation*}

\begin{equation*}
~~4 \equiv \left(
\begin{array}{ccc}
1 & 1 & 1 \\
0 & 1 & 0 \\
0 & 0 & 1 \\
\end{array}
\right),~~
5 \equiv \left(
\begin{array}{ccc}
0 & 1 & 1 \\
0 & 0 & 0 \\
0 & 0 & 0 \\
\end{array}
\right),~~
6 \equiv \left(
\begin{array}{ccc}
0 & 0 & 1 \\
0 & 0 & 0 \\
0 & 0 & 0 \\
\end{array}
\right),~~
7 \equiv \left(
\begin{array}{ccc}
1 & 0 & 1 \\
0 & 1 & 0 \\
0 & 0 & 1 \\
\end{array}
\right),
\end{equation*}

\begin{equation*}
~~8 \equiv \left(
\begin{array}{ccc}
0 & 0 & 0 \\
0 & 1 & 0 \\
0 & 0 & 1 \\
\end{array}
\right),~~
9 \equiv \left(
\begin{array}{ccc}
1 & 0 & 0 \\
0 & 0 & 0 \\
0 & 0 & 0 \\
\end{array}
\right),~
10 \equiv \left(
\begin{array}{ccc}
1 & 1 & 0 \\
0 & 0 & 0 \\
0 & 0 & 0 \\
\end{array}
\right),~
11 \equiv \left(
\begin{array}{ccc}
0 & 1 & 0 \\
0 & 1 & 0 \\
0 & 0 & 1 \\
\end{array}
\right),
\end{equation*}

\begin{equation*}
~12 \equiv \left(
\begin{array}{ccc}
1 & 1 & 1 \\
0 & 0 & 0 \\
0 & 0 & 0 \\
\end{array}
\right),~
13 \equiv \left(
\begin{array}{ccc}
0 & 1 & 1 \\
0 & 1 & 0 \\
0 & 0 & 1 \\
\end{array}
\right),~
14 \equiv \left(
\begin{array}{ccc}
0 & 0 & 1 \\
0 & 1 & 0 \\
0 & 0 & 1 \\
\end{array}
\right),~
15 \equiv \left(
\begin{array}{ccc}
1 & 0 & 1 \\
0 & 0 & 0 \\
0 & 0 & 0 \\
\end{array}
\right),
\end{equation*}
we see that 0 is the additive identity, 1 is the multiplicative identity and the only invertible elements  are 1, 2, 4 and 7. 
The ring features two (two-sided) maximal ideals, namely

\begin{equation*}
I_l = \{0,3,5,6,8,11,13,14\}
\end{equation*}
and
\begin{equation*}
I_r = \{0,3,5,6,9,10,12,15\},
\end{equation*}
and its Jacobson radical reads
\begin{equation*}
J = \{0,3,5,6\}.
\end{equation*}
From the above-given matrix representation of $R$ we find that $R^2$ contains nine distinct free cyclic submodules generated by nonunimodular vectors, which are listed in Table 1.

\begin{table}[pth!]
\begin{center}
\caption{The explicit form of the nine nonunimodular free cyclic submodules. } \vspace*{0.4cm}
\resizebox{\columnwidth}{!}{%
\begin{tabular}{|r|r|r|r|r|r|r|r|r|r|}
\hline \hline
  $\alpha$ & $R(3,8)$ & $R(5,8)$ & $R(6,8)$ & $R(8,11)$  & $R(8,13)$  & $R(8,14)$  & $R(8,6)$  & $R(8,5)$  & $R(8,3)$   \\
\hline
  0        & $(0,0)$  & $(0,0)$  & $(0,0)$  &  $(0,0)$   & $(0,0)$    & $(0,0)$    & $(0,0)$   & $(0,0)$   & $(0,0)$   \\
	1        & $(3,8)$  & $(5,8)$  & $(6,8)$  &  $(8,11)$  & $(8,13)$   & $(8,14)$   & $(8,6)$   & $(8,5)$   & $(8,3)$   \\
	2        & $(3,11)$ & $(5,11)$ & $(6,11)$ &  $(11,8)$  & $(11,14)$  & $(11,13)$  & $(11,6)$  & $(11,5)$  & $(11,3)$   \\
  3        & $(0,3)$  & $(0,3)$  & $(0,3)$  &  $(3,3)$   & $(3,3)$    & $(3,3)$    & $(3,0)$   & $(3,0)$   & $(3,0)$   \\
  4        & $(3,13)$ & $(5,13)$ & $(6,13)$ &  $(13,14)$ & $(13,8)$   & $(13,11)$  & $(13,6)$  & $(13,5)$  & $(13,3)$   \\
  5        & $(0,5)$  & $(0,5)$  & $(0,5)$  &  $(5,5)$   & $(5,5)$    & $(5,5)$    & $(5,0)$   & $(5,0)$   & $(5,0)$   \\
  6        & $(0,6)$  & $(0,6)$  & $(0,6)$  &  $(6,6)$   & $(6,6)$    & $(6,6)$    & $(6,0)$   & $(6,0)$   & $(6,0)$   \\
  7        & $(3,14)$ & $(5,14)$ & $(6,14)$ &  $(14,13)$ & $(14,11)$  & $(14,8)$   & $(14,6)$  & $(14,5)$  & $(14,3)$   \\
  8        & $(0,8)$  & $(0,8)$  & $(0,8)$  &  $(8,8)$   & $(8,8)$    & $(8,8)$    & $(8,0)$   & $(8,0)$   & $(8,0)$   \\
  9        & $(3,0)$  & $(5,0)$  & $(6,0)$  &  $(0,3)$   & $(0,5)$    & $(0,6)$    & $(0,6)$   & $(0,5)$   & $(0,3)$   \\
 10        & $(3,3)$  & $(5,3)$  & $(6,3)$  &  $(3,0)$   & $(3,6)$    & $(3,5)$    & $(3,6)$   & $(3,5)$   & $(3,3)$   \\
 11        & $(0,11)$ & $(0,11)$ & $(0,11)$ &  $11,11)$  & $(11,11)$  & $(11,11)$  & $(11,0)$  & $(11,0)$  & $11,0)$   \\
 12        & $(3,5)$  & $(5,5)$  & $(6,5)$  &  $(5,6)$   & $(5,0)$    & $(5,3)$    & $(5,6)$   & $(5,5)$   & $(5,3)$   \\
 13        & $(0,13)$ & $(0,13)$ & $(0,13)$ &  $(13,13)$ & $(13,13)$  & $(13,13)$  & $(13,0)$  & $(13,0)$  & $(13,0)$   \\
 14        & $(0,14)$ & $(0,14)$ & $(0,14)$ &  $(14,14)$ & $(14,14)$  & $(14,14)$  & $(14,0)$  & $(14,0)$  & $(14,0)$   \\
 15        & $(3,6)$  & $(5,6)$  & $(6,6)$  &  $(6,5)$   & $(6,3)$    & $(6,0)$    & $(6,6)$   & $(6,5)$   & $(6,3)$   \\
 \hline \hline
\end{tabular}}%
\end{center}
\end{table}
We shall show that the way how these free cyclic submodules are interwoven is intricately related to the structure of the generalized quadrangle of order two, the doily. To this end we employ a duad-syntheme model of the latter (see, for example, \cite{paythas}). Take a six-element set  $S=\{1,2,3,4,5,6\}$. Let us call a two-element subset of $S$ a duad and a set of three duads forming a partition of $S$ a syntheme. Then the point-line incidence structure whose points are 15 duads and whose lines are 15 synthemes, with incidence being containment, is isomorphic to the doily. The structure of the doily is illustrated in Figure 1 -- {\it left}. Here, the points of the doily are represented by circles, labeled by duads, and its lines are represented by nine straight segments, three concentric circles and three arcs of circles; one can readily check that each line corresponds to a syntheme. Next, we employ the following bijection between the 15 duads and 15 nontrivial vectors $(a,b) \in R^2$, where $a,b \in J$ and $(a,b) \neq (0,0)$: 

\begin{equation*}
\begin{split}
\{1,2\} \leftrightarrow (3,3),~~\{1,3\} \leftrightarrow (5,3),~~\{1,4\} \leftrightarrow (0,6),~~\{1,5\} \leftrightarrow (3,6),~~\{1,6\} \leftrightarrow (5,0),\\ 
\{2,3\} \leftrightarrow (6,0),~~\{2,4\} \leftrightarrow (3,5),~~\{2,5\} \leftrightarrow (0,5),~~\{2,6\} \leftrightarrow (6,3),~~\{3,4\} \leftrightarrow (5,5),\\
\{3,5\} \leftrightarrow (6,5),~~\{3,6\} \leftrightarrow (0,3),~~\{4,5\} \leftrightarrow (3,0),~~\{4,6\} \leftrightarrow (5,6),~~\{5,6\} \leftrightarrow (6,6).
\end{split}
\end{equation*}
\noindent
We thus get a new labeling of the points of the doily in terms of these particular nonunimodular vectors of $R^2$, as illustrated in Figure 1 -- {\it right}. From comparison of the latter figure with Table 1 it follows that each submodule shares with the doily seven vectors forming three concurrent lines -- as depicted in Figure 2. From this figure one can easily discern that six lines of the doily have a different 
status than the other nine, as each of them belongs to three submodules. Given the fact that each point of concurrence belongs to two such lines, the nine points and the six lines are found to form inside the doily a point-line incidence structure isomorphic to the generalized quadrangle of type GQ(2,1) (see \cite{paythas}).

\begin{figure}[t]
\centering
\includegraphics[width=4cm]{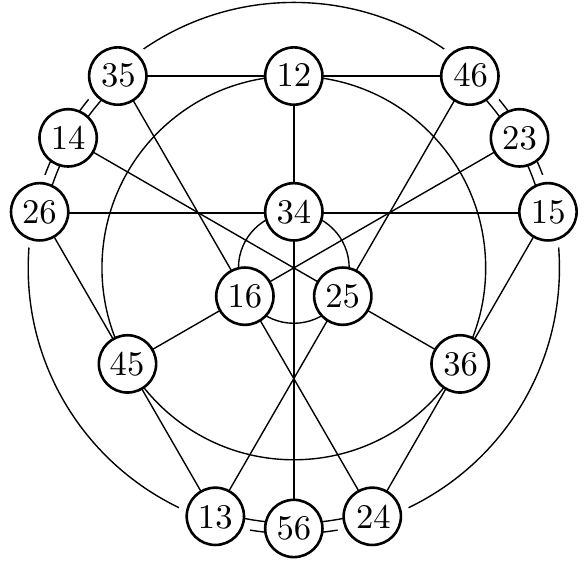} \hspace*{3.5cm}
\includegraphics[width=4cm]{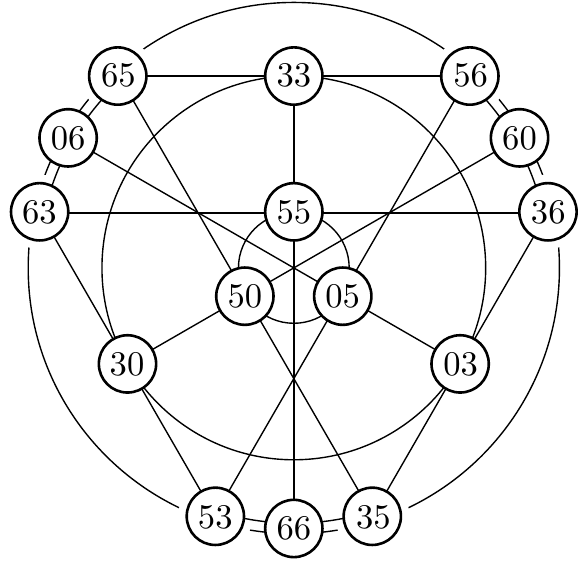}
\caption{A pictorial representation of the doily, with its points labeled by duads ({\it left}) and nonunimodular vectors ({\it right}); both a duad $\{a,b\}$ and a vector $(a,b)$ are abbreviated to $ab$.}
\end{figure}

\begin{figure}[pth!]
\centering
\includegraphics[width=11cm]{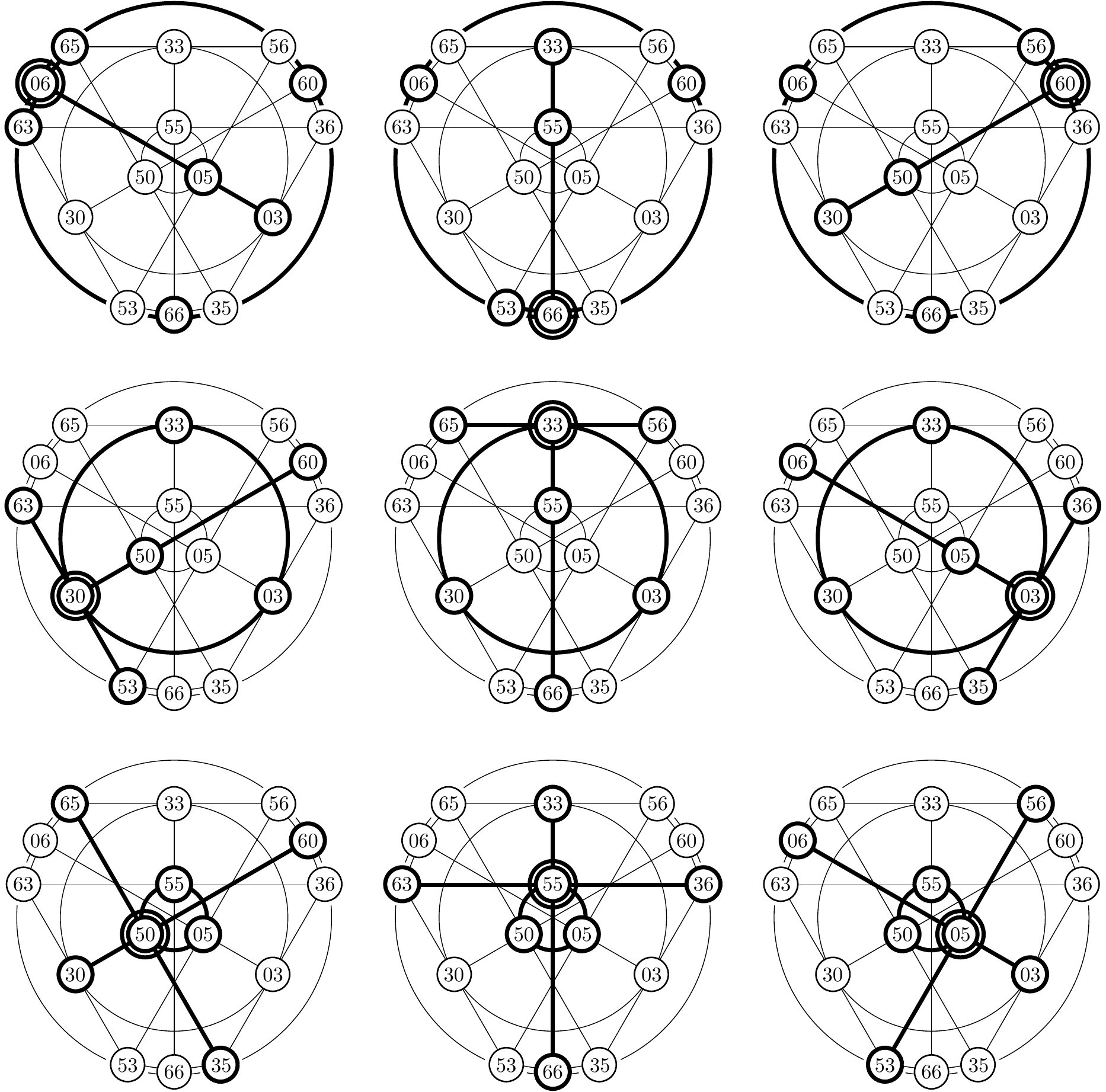}
\caption{A graphical illustration of `Jacobson traces' of individual nonunimodular free cyclic submodules in the core doily; consecutively from top left to bottom right there are shown the traces of $R(6,8)$, $R(8,14)$, $R(8,6)$, $R(8,3)$, $R(8,11)$, $R(3,8)$, $R(8,5)$, $R(8,13)$ and $R(5,8)$. In each subfigure the corresponding concurrent lines are shown in boldface, the point of concurrence being encircled.}
\end{figure}

A complete view of the relation between individual submodules is outlined in Figure 3. The pronounced automorphism of order three of the figure stems from the fact that the submodules form three disjoint triples according to the number of shared vectors. One further observes that all vectors lying on our submodules acquire values from the ideal $I_l$. It can readily be verified that a completely analogous geometric structure is obtained if we take the free {\it right} module, in which case the corresponding vectors have entries from the ideal $I_r$. Obviously, the two structures share the same doily, as its points are labeled by vectors from $J^2$. At this point it is well worth recalling the existence of a similar geometrical structure in the case of the smallest ring of ternions and its three-dimensional free left (and also right) module \cite{sanetal}.
There the associated geometry, referred to as the `Fano-snowflake,' has its center isomorphic to the Fano plane, a generalized triangle of order two (see also \cite{hasa} for generalization to an arbitrary ring of ternions).

\begin{figure}[pth!]
\centering
\includegraphics[width=14cm]{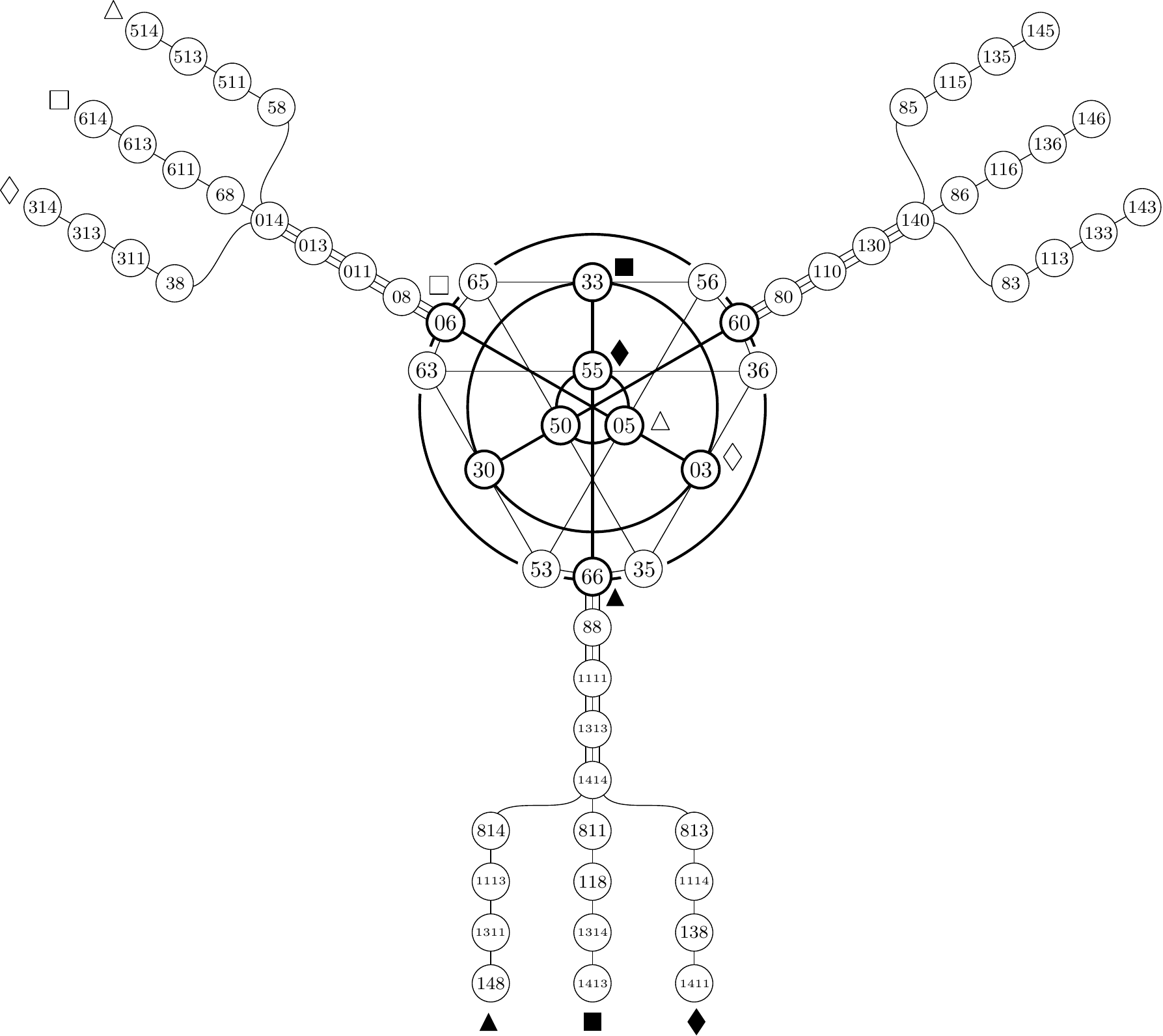}
\caption{Visualisation of the full geometric structure formed by the nine nonunimodular free cyclic submodules, with the doily lying in its center. For each submodule there are shown all of its vectors except for the trivial one. For six submodules pairs of identical symbols are employed to identify the corresponding points of concurrence; the remaining three cases are obtained by swapping the figure with respect to the vertical axis passing through the center. The distinguished GQ(2,1) of the doily is shown in bold.}
\end{figure}

\pagebreak
The occurrence of the doily in this remarkable nonunimodular ring-theoretic setting is quite intriguing also in view of possible physical applications. For example, among the finite geometric concepts relevant for the theory of quantum information, the doily -- though in various disguises -- has been recognized to play the foremost role. First, being isomorphic to the symplectic polar space of type  $W(3,2)$, it underlies the commutation relations between the elements of the two-qubit Pauli group \cite{twoqub} and provides us with simplest settings (namely GQ(2,\,1)'s) for observable proofs of quantum contextuality.  Second, being isomorphic to a non-singular quadric of type $\mathcal{Q}(4,2)$, it also lies in the heart of a remarkable magic three-qubit Veldkamp line of form theories of gravity and its four-qubit extensions \cite{mvl}. Finally, being a subquadrangle of a generalized quadrangle of type GQ(2,\,4), it enters in an essential way certain black-hole entropy formulas and the so-called black-hole/qubit correspondence \cite{bhq}.
We, therefore, believe that the above-described doily-based geometry as a whole will eventually find its way into (quantum) physics as well.

\section*{Acknowledgments}
This work was supported, in part, by the Slovak VEGA Grant Agency, Project $\#$ 2/0003/16, and by the National Scholarship Programme of the Slovak Republic. We are extremely grateful to Zsolt Szab\'o for electronic versions of the figures.

\end{document}